\documentclass{amsart}

\usepackage{amsmath,amsthm,amssymb,amscd}
\usepackage{amsfonts}
\usepackage{rotating}
\usepackage{euscript}
\usepackage{pst-node}
\usepackage{epsfig,verbatim}
\def\be{\begin{equation}}
\def\ee{\end{equation}}

\def\C{{\mathbb C}} 
\def\f{\EuScript}
 
\def\P{{\mathbb P}}

\def\e{\eqref}
\def\phi{{\varphi}}
 
\def\tt{\widetilde}
\def\deg{{\rm deg\,}}

\def\GCD{{\rm GCD }}

\def\bp{\begin{proposition}}
\def\ep{\end{proposition}}

\def\bt{\begin{theorem}}
\def\et{\end{theorem}}
\def\br{\begin{remark}}
\def\er{\end{remark}}
\def\be{\begin{equation}}
\def\bee{\begin{equation*}}
\def\l{\label}

\def\e{\eqref}

\def\ee{\end{equation}}
\def\eee{\end{equation*}}
\def\bl{\begin{lemma}}
\def\el{\end{lemma}}
\def\bc{\begin{corollary}}
\def\ec{\end{corollary}}
\def\pr{\noindent{\it Proof. }}

\def\bd{\begin{definition}}
\def\ed{\end{definition}}
\def\t{\widetilde}

\newtheorem{theorem}{Theorem}[section]
\newtheorem{lemma}{Lemma}[section]
\newtheorem{definition}{Definition}[section]
\newtheorem{corollary}{Corollary}[section]
\newtheorem{proposition}{Proposition}[section]
\newtheorem{remark}{Remark}[section]

\usepackage{mathtools}

\mathtoolsset{showonlyrefs}

\begin{document}
\title{Semiconjugate rational functions: a dynamical approach}
\author{F. Pakovich}

\begin{abstract}
Using dynamical methods  we give a new proof of the theorem saying that if $A,B,X$ are rational functions of variable $z$ of degree at least two  such that $A\circ X=X\circ B$ and  $\C(B,X)=\C(z)$, then 
the Galois closure of the field extension  $\C(z)/\C(X)$ has genus zero or one.  
\end{abstract} 
\maketitle

\begin{section}{Introduction}
Let $A$ and $B$ be  rational functions of degree at least two on the Riemann sphere. 
The function $B$ is said to be semiconjugate to the function $A$
if there exists a non-constant rational function $X$
such that \be \l{i1} A\circ X=X\circ B.\ee
Notice that for $\deg X=1$  condition \eqref{i1} reduces to the usual conjugacy condition
while for $B=A$ it reduces to the commutativity condition  
\be \l{i2}
A\circ X=X\circ A.
\ee

A solution of equation \eqref{i1} is called 
primitive if the functions $X$ and $B$ generate the whole field of rational functions $\C(z)$. Up to a certain degree, 
the description of solutions of \eqref{i1}  reduces to the description of primitive solutions. 
Indeed, 
by the L\"uroth theorem, the field 
$\C(X,B)$ is generated by some rational function $W$.  
Therefore, if  $\C(X,B)\neq \C(z)$, then 
 there exists a rational function $W$ of degree greater than one such that
\be \l{of} B=\tt B\circ W, \ \ \ X=\tt X\circ W\ee 
for some rational functions $\tt X$ and $\tt B$. Substituting now \eqref{of} in \eqref{i1} we see  that the triple $A, \tt X,W\circ \tt B$ is another solution of \e{i1}. This new solution is not necessary primitive, however 
$\deg \tt X<\deg X$. Therefore, after a finite number of similar transformations we will arrive to a primitive solution.

Semiconjugate  rational functions were investigated  at length in the series of papers \cite{semi}, \cite{rec},  \cite{fin}, \cite{lattes}. In particular, it was shown in \cite{semi} that all primitive solutions of 
\eqref{i1} are related to discrete automorphism groups of $\C$ and $\C\P^1$,
implying that  corresponding 
functions $X$  have a very restricted form. 
Recall that for a rational function $X$  
its  normalization $\t X$ is defined as a holomorphic function of the lowest possible degree
between compact Riemann surfaces  $\t X:\,\t S_X\rightarrow \C\P^1$  such that $\t X$ is a Galois covering and
 $$\t X=X\circ H$$ for some  holomorphic map $H:\,\t S_X\rightarrow \C\P^1$. From the algebraic point of view the passage from $X$ to $\t X$ 
corresponds to the passage from the field extension $\C(z)/\C(X)$ to its Galois closure. 
In these terms, the main result of \cite{semi} about primitive solutions of \eqref{i1} may be formulated as follows.

\bt \l{main} Let $A,B,X$ be rational functions of degree at least two  such that $A\circ X=X\circ B$ and  $\C(B,X)=\C(z)$.   
Then the Galois closure of the field extension $\C(z)/\C(X)$ has genus zero or one.    
\et 

Observe a similarity between this result and the Ritt theorem (\cite{r}) saying that if
two rational functions $A$ and $X$ commute and have no iteration in common, then $A$ and $X$ 
either are  Latt\`es functions, or are conjugate to powers or Chebyshev polynomials. Indeed,  powers and Chebyshev polynomials are the simplest examples of rational functions such that 
$g(\t S_X)=0$. On the other hand, Latt\`es maps are examples of rational functions with $g(\t S_X)=1$. 
Rational functions $X$ with $g(\t S_X)=0$ can be listed explicitly, while functions with $g(\t S_X)=1$
admit a simple geometric description (see \cite{gen}). 
Notice that rational functions with $g(\t S_X)\le 1$
can be described 
through their ramification,
implying that Theorem \ref{main} is  equivalent to Theorem 6.1 of \cite{semi} (see Section 5 below).

The problem of describing commuting and semiconjugate rational functions naturally belongs to  dynamics (see e.g. the  papers \cite{be},  \cite{e}, \cite{e2}, \cite{f}, \cite{j}, \cite{ms}, \cite{pj}). 
In particular, in the papers of Fatou and Julia \cite{f}, \cite{j} commuting rational functions  were investigated by dynamical methods, 
requiring however an assumption that the
Julia sets of considered functions do not coincide with the whole Riemann sphere. On the other hand, the Ritt theorem about commuting rational functions cited above 
was proved by non-dynamical methods. 
In his paper, Ritt remarked that ``it would be interesting to know whether a proof can also be effected by the use of Poncar\'e functions employed by Julia''.
Sixty six years later such a proof was given   
by Eremenko  \cite{e2}. Notice that the Ritt theorem also follows from the results of \cite{fin} 
about solutions of equation \eqref{i1} with fixed $B$.

Similarly to the paper of Ritt \cite{r}, the paper  \cite{semi}  does not use any dynamical methods, but relies
on a  study of maps between two-dimensional orbifolds 
associated with rational functions.
At the same time, it is interesting to find approaches to equation \eqref{i1}
 involving ideas from  dynamics, 
and the goal of this paper is to provide a ``dynamical'' proof of Theorem \ref{main}. 
In fact, we give {\it two} such proofs. The first one  exploits a link between equation \eqref{i1} and  Poincar\'e functions.
The second one is based on the interpretation of $\t S_X$ as an invariant curve for the dynamical system \be \l{ds} (x_1,x_2,\dots, x_n)\rightarrow (A(x_1),A(x_2),\dots ,A(x_n))\ee on $(\C\P^1)^n$.
The last proof is inspired by the recent paper \cite{ms} describing invariant varieties for dynamical systems of the form 
  $$(x_1,x_2,\dots, x_n)\rightarrow (C_1(x_1),C_2(x_2),\dots ,C_n(x_n)),$$ where $C_1,C_2,\dots, C_n$ are {\it polynomials}, and relating such varieties with polynomial solutions of \eqref{i1}. The analysis of  equation \eqref{i1} in the paper \cite{ms}, based on the Ritt theory of polynomial decompositions (\cite{r1}),   
 does not extend to arbitrary rational functions. Nevertheless, the relation between the semiconjugacy condition and invariant varieties established in \cite{ms} suggests that there should be some interpretation of the results of  \cite{semi} in terms of	dynamical systems of form  \eqref{ds}, and we show that this is indeed the case.

The paper is organized as follows. In the second section we recall the description of  $\t S_X$ in terms  of algebraic equations, and 
give a criterion for a rational function $X$ to satisfy the condition $g(\t S_X)\leq 1$.
In the third and the fourth sections we provide two proofs of Theorem 1.1 using two approaches described above.  
Finally, in the fifth section we show that Theorem \ref{main} is equivalent to  Theorem 6.1 of \cite{semi} which describes primitive solutions of \eqref{i1} in terms of orbifolds.

\end{section}

\begin{section}{Meromorphic parametrizations and normalizations} 

Let $\f C$ be an  irreducible algebraic curve in $\C^n$.
Recall that a meromorphic parametrization of $\f C$ on $\C$ is a collection of functions  $\psi_1, \psi_2, \dots ,\psi_n$ such that  

\begin{itemize} 
\vskip 0.1cm
\item  $\psi_1,$ $\psi_2,$  $\dots,$ $\psi_n$ are non-constant and meromorphic on $\C$, \vskip 0.1cm
\item  $(\psi_1(z), \psi_2(z), \dots ,\psi_n(z))\in \f C$ whenever $\psi_i(z)\neq \infty$, $1\leq i \leq n,$ 
\vskip 0.1cm
\item with finitely many exceptions, every point of $\f C$ is of the form  \linebreak  $(\psi_1(z), \psi_2(z), \dots ,\psi_n(z))$ for some $z\in \C.$ 
\end{itemize} 
Notice that the last condition in fact is a corollary of the first two. 

By the classical theorem of Picard (\cite{p}), a plane algebraic curve $\f C$ which can be parametrized by functions meromorphic on $\C$  has genus zero or one  (see  e.g. \cite{bn}). 
We will use the following slightly more general version of this theorem which can proved in the same way. 

\bt \l{picard} If  an irreducible algebraic curve  $\f C$ in $\C^n$ has a meromorphic para\-met\-rization on $\C$, then  $\f C$ has genus zero or one. \qed
\et

 Let $X:\C\P^1\rightarrow \C\P^1$ be a rational function of degree $d$.  The normalization $\t X: \t S_X \rightarrow \C\P^1$ 
can be described by the following construction (see \cite{fried}, $\S$I.G).
Consider the fiber product of the cover $X:\C\P^1\rightarrow \C\P^1$ with itself $d$ times, that is a subset $\f L$ of $(\C\P^1)^d$ consisting of $d$-tuples  with a common image 
under $X$. Clearly, $\f L$ is an algebraic variety of dimension one 
 defined by the algebraic equations \be \l{uri} X(z_i)-X(z_j)=0, \ \ \ \ 1\leq i,j\leq d,\ \ \  i\neq j.\ee
Let ${\f L}_0$ be a variety obtained from $\f L$ by removing the components where two or more coordinates coincide,  $\f N$   an irreducible component of 
${\f L}_0$, and ${\f N}^{\prime}\xrightarrow{\pi^{\prime}} \f N$ the desingularization map. In this notation the following statement holds.

\bt \l{fr} 
The map $\psi: {\f N^{\prime}}\rightarrow \C\P^1$ given by the composition
\be \l{pro} {\f N^{\prime}}\xrightarrow{\pi^{\prime}} {\f N}\xrightarrow{\pi_i}\C\P^1\xrightarrow{X}\C\P^1,\ee
where ${\f N}$ is 
any irreducible component of ${\f L}_0$ and $\pi_i$ is the projection to any coordinate, is the normalization of $X$. \qed
\et

Combining Theorem \ref{picard} and Theorem \ref{fr} we obtain the following characterization of rational functions $X$ with   $g(\t S_X)\leq 1$.
 
\bt \l{param} Let $X$ be a rational function of degree $d$. Then $g(\t S_X)\leq 1$ if and only if there exist $d$ distinct functions $\psi_1, \psi_2, \dots ,\psi_d$ meromorphic on $\C$ such that 
\be \l{u} X(\psi_i)-X(\psi_j)=0,  \ \ \ \ 1\leq i,j\leq d,\ \ \  i\neq j.\ee  
\et
\pr Equalities \eqref{u} imply that some irreducible component $\f N$ of $\f L_0$ admits a meromorphic parametrization.
Since ${\f N^{\prime}}=\t S_X$ by Theorem \ref{fr}, it follows  from Theorem \ref{picard}
that $g(\t S_X)\leq 1$.  

In the other direction, if  $g(\t S_X)\leq 1$, then  
taking different coordinate projections in \eqref{pro}  
we obtain $d$ distinct functions
$$\theta_i=\pi_i\circ \pi^{\prime},\ \ \ \  1\leq i \leq d,$$
from $\t S_X$ to $\C\P^1$ 
such that 
$$X(\theta_i)-X(\theta_j)=0,  \ \ \ \ 1\leq i,j\leq d,\ \ \  i\neq j.$$ If $g(\t S_X)=0,$ these functions are rational and therefore meromorphic on $\C.$ 
On the other hand, if $g(\t S_X)=1,$  we obtain meromorphic functions satisfying \eqref{u}
setting 
$$\psi_i= \theta_i\circ\tau, \ \ \ \ 1\leq i \leq d,$$ where $\tau:\C \rightarrow \t S_X$ is the universal covering of $\t S_X$.
 \qed

\end{section}

\begin{section}{Semicoconjugate functions and Poincar\'e functions
}

Let $A$ be a rational function and $z_0$ its repelling fixed point. 
Recall that the Poincar\'e function $\f P_{A,z_0}$ associated with $z_0$ is a function meromorphic on $\C$ such that 
$\f P_{A,z_0}(0)=z_0$,  $\f P_{A,z_0}^{\prime}(0)=1$, and the diagram    
\be \l{xor}
\begin{CD} 
\C @>\lambda z>> \C \\
@VV \f P_{A,z_0} V @VV \f P_{A,z_0} V\\ 
\f \C\P^1 @>A>> \f\C\P^1\
\end{CD}
\ee
commutes. The 
Poincar\'e function always exists and is unique  (see e.g \cite{mil}).

\bl \l{l} Let $X$ and $B$ be rational functions such that $\C(X,B)=\C(z).$ Then  for all but finitely many $z_0\in \C$ 
the set $B(X^{-1}\{z_0\})$ contains $\deg X$ distinct points.
\el 
\pr Since $\C(X,B)=\C(z)$, there exist $U,V\in \C[x,y]$ such that  
$$z=\frac{U(X,B)}{V(X,B)}.$$ This implies that for $z_1\neq z_2$ such that $X(z_1)=X(z_2)$ the inequality 
$B(z_1)\neq B(z_2)$ holds, unless $z_1$ or $z_2$ is a zero of the polynomial  $V(X,B).$  
Therefore, if $z_0$ is neither a critical value of $X$ nor an $X$-image of a zero of  $V(X,B)$, 
the set $B(X^{-1}\{z_0\})$ contains $\deg X$ distinct points, since 
 $X^{-1}\{z_0\}$ contains $\deg X$ distinct points and their $B$-images are distinct.
\qed

\vskip 0.2cm

Combining the uniqueness of the Poincar\'e function with Theorem \ref{param} we can prove 
Theorem \ref{main} as follows. Let $A,$ $B,$ and $X$ be rational functions of degree at least two such that 
 the diagram    
\be \l{i11}
\begin{CD} 
\C\P^1 @> B>> \C\P^1 \\
@VV X V @VV X V\\ 
\f \C\P^1 @>A>> \f\C\P^1\
\end{CD}
\ee commutes and  $\C(X,B)=\C(z)$. 
Since the number of repelling periodic points of $A$ is infinite, it follows from Lemma \ref{l} that 
we can find a repelling periodic point $z_0\in \C$  such that for any point $z$ in the forward $A$-orbit of $z_0$ the set $B(X^{-1}\{z\})$ contains $\deg X$ distinct points.
Since \eqref{i11}
implies that $$B(X^{-1}\{z_0\})\subseteq X^{-1}\{A(z_0)\},$$ this yields  that
 $$B(X^{-1}\{z_0\})= X^{-1}\{A(z_0)\},$$ and,  
 inductively, that 
  $$B^{\circ k}(X^{-1}\{z_0\})= X^{-1}\{A^{\circ k}(z_0)\}, \ \ \ \ k\geq 1.$$ 
 In particular, for $n$ equal to the period of $z_0$ we have:
$$B^{\circ n}(X^{-1}\{z_0\})=X^{-1}\{z_0\}.$$
Therefore, the restriction 
of the rational function $B^{\circ n}$ on the set $X^{-1}\{z_0\}$ is a permutation of its elements, and hence  
for certain $l\geq 1$ all the points of $X^{-1}\{z_0\}$ are fixed points of $B^{\circ nl}$. Thus, considering instead of $A$ and $B$ their iterations we can assume that $z_0$ is a fixed point of $A$ and the set $X^{-1}\{z_0\}$
consists of $d=\deg X$ distinct fixed points $z_1,z_2,\dots , z_d$ of $B.$ 

Since the points $z_1,z_2,\dots , z_d$ are not critical points of $X$, the map $X$ is invertible near each of them implying that the multipliers of $B$ at  $z_1,z_2,\dots , z_d$ are all equal 
to the multiplier $\lambda$ of $A$ at $z_0$, so that   $z_1,z_2,\dots , z_d$ are repelling fixed points of $B$. Clearly, for each $i,$ $1\leq i \leq d,$ 
we can complete commutative diagram \eqref{i11} to  the commutative diagram 
$$
\begin{CD}
\C @>\lambda z>> \C \\
@VV \f P_{B,z_i} V @VV \f P_{B,z_i} V\\ 
\f \C\P^1 @>B>> \f\C\P^1\\
@VV X V @VV X V\\ 
\C\P^1 @> A >> \C\P^1 \ ,
\end{CD}
$$ where $\f P_{B,z_i}$, $1\leq i \leq d,$  is  the corresponding Poincar\'e function for $B$. 
Since the functions $X\circ \f P_{B,z_i},$ $1\leq i \leq d,$ are meromorphic, 
it follows now from the uniqueness of the Poincar\'e function that there exist $\alpha_1,\alpha_2,\dots\alpha_n\in \C\setminus\{0\}$ such that 
\be \l{a} \f P_{A,z_0}(z)=X\circ \f P_{B,z_1}(\alpha_1z)=X\circ \f P_{B,z_2}(\alpha_2z)=\dots = X\circ \f P_{B,z_d}(\alpha_dz).\ee 
Moreover, the functions  $\f P_{B,z_i}(\alpha_iz),$ $1\leq i \leq d,$ are distinct
 since the points $$z_i=\f P_{B,z_i}(0), \ \ \ 1\leq i \leq d,$$ 
are distinct. 
Applying now Theorem \ref{param} to  equality \eqref{a}, we see that  $g(\t S_X)\leq 1$.

\end{section}

\begin{section}{Semicoconjugate functions and invariant curves
}

Let  $R_1,R_2,\dots , R_d$ be rational functions, and let $\f R:(\C\P^1)^d\rightarrow (\C\P^1)^d$ be the map
$$(x_1,x_2,\dots, x_d)\rightarrow (R_1(x_1),R_2(x_2),\dots ,R_d(x_d)).$$ Say that an algebraic curve $\f C$ in $(\C\P^1)^d$ is $\f R$-invariant if 
$\f R(\f C)=\f C.$ Invariant curves possess the following property (cf. \cite{ms}, Proposition 2.34). 

\bt \l{ge} Let $R_1,R_2,\dots , R_d$  be rational functions of degree at least two and $\f C$ an  irreducible  $\f R$-invariant curve in $(\C\P^1)^d.$ Then $g(\f C)\leq 1.$ 
\et 
\pr Since $\f C$ is $\f R$-invariant, the map $\f R$ lifts to a holomorphic map 
$${\f R}^{\prime}:\, {\f C}^{\prime}\rightarrow {\f C}^{\prime},$$ where $\f C^{\prime}$ is a desingularization of $\f C.$
Applying now 
the Riemann-Hurwitz formula
$$2g({\f C}^{\prime})-2=(2g({\f C}^{\prime})-2)\deg {\f {\f R}^{\prime}}+\sum_{P\in \f C^{\prime}}(e_p-1),$$ we see  that 
$g(\f C^{\prime})\leq 1$, unless $\deg {\f R}^{\prime}=1.$ 

Furthermore, for $\deg {\f R}^{\prime}=1$ the inequality $g({\f C}^{\prime})\leq 1 $ still holds. Indeed, since the automorphism group of a Riemann surface of genus greater than one is finite, if $g({\f C}^{\prime})\geq 2 $, then for some $k\geq 1$ the map
$({\f R}^{\prime})^{\circ k}$ is the identical automorphism of $\f C^{\prime}$, implying that the maps  
$$(z_1,z_2,\dots, z_d)\rightarrow R_i^{\circ k}(z_i), \ \ \ \ \ 1\leq i \leq d,$$ are identical on $\f C$. 
However, since each $R_i,$ $1\leq i \leq d,$ has degree at least two, in this case for every point  of $\f C$ its 
$i$th coordinate belongs to a finite subset of $\C\P^1$ consisting of   
fixed point of $R_i^{\circ k}$, implying that $\f C$ is a finite set. \qed

\vskip 0.2cm

Using Theorem \ref{ge} and Theorem \ref{fr} we obtain a proof of Theorem \ref{main} as follows. 
Define the maps $\f A$, $\f B$, and $\f X$ from  $(\C\P^1)^d$ to  $(\C\P^1)^d$ by the formulas 
$$\f A:\, (x_1,x_2,\dots, x_d)\rightarrow (A(x_1),A(x_2),\dots ,A(x_d)),$$ 
$$\f B:\, (x_1,x_2,\dots, x_d)\rightarrow (B(x_1),B(x_2),\dots ,B(x_d)),$$
$$\f X:\, (x_1,x_2,\dots, x_d)\rightarrow (X(x_1),X(x_2),\dots ,X(x_d)).$$
Clearly, equality \eqref{i1} implies that the diagram 
\be \l{kaba}
\begin{CD}
(\C\P^1)^d @>\f B>> (\C\P^1)^d \\
@VV\f X V @VV\f X V\\ 
(\C\P^1)^d @>\f A >> (\C \P^1)^d
\end{CD}
\ee
commutes. By construction, the variety $\f L$ defined by equations \eqref{uri} is the preimage of the diagonal $\Delta$ in $(\C\P^1)^d$ 
under the map  $\f X:(\C\P^1)^d\rightarrow(\C\P^1)^d$.
 Therefore, since $\f A(\Delta)=\Delta$, it follows from \eqref{kaba} that $\f B(\f L)\subseteq\f L$. 
Moreover, Lemma \ref{l} implies that $\f B(\f L_0)\subseteq \f L_0$. 
Since $\f L_0$ has a finite number of irreducible components, this implies that there exists  an irreducible component $\f N_0$ of $\f L_0$ such that  $\f B^{\circ k}(\f N_0)= \f N_0$ for some $k\geq 1.$
Since by Theorem \ref{fr} the equality  $g(\f N_0)=g(\t{S}_X)$ holds,  it follows now from Theorem \ref{ge} that $g( \t{S}_X)\leq 1.$

\end{section}

\begin{section}{Semicoconjugate functions and orbifolds} 
Recall that an orbifold $\f O$ on $\C\P^1$ is 
a ramification function $\nu:\C\P^1\rightarrow \mathbb N$ which takes the value $\nu(z)=1$ except at finite number of points. The Euler characteristic of an orbifold $\f O$ is defined by the formula  
\be \l{char}  \chi(\f O)=2+\sum_{z\in \C\P^1}\left(\frac{1}{\nu(z)}-1\right).\ee
A rational function $f$ is called a covering map $f:\,  \f O_1\rightarrow \f O_2$
between orbifolds
$\f O_1$ and $\f O_2$
if for any $z\in \C\P^1$ the equality 
\be \l{us} \nu_{2}(f(z))=\nu_{1}(z)\deg_zf\ee holds, where $\deg_zf$ denotes the local degree of $f$ at the point $z$. 
If $f:\,  \f O_1\rightarrow \f O_2$ is a covering map, then  
the Riemann-Hurwitz formula implies that \be\l{rh} \chi(\f O_1)=\deg  f\cdot\chi(\f O_2).\ee 
In case if a weaker than \eqref{us}  condition
\be \l{rys} \nu_{2}(f(z))=\nu_{1}(z)\GCD(\deg_zf, \nu_{2}(f(z))\ee 
holds,  $f$ is called a minimal holomorphic map 
between orbifolds $\f O_1$ and $\f O_2$. 

With each rational function $f$ 
one can associate in a natural way two orbifolds $\f O_1^f$ and 
$\f O_2^f$, setting $\nu_2^f(z)$  
equal to the least common multiple of the local degrees of $f$ at the points 
of the preimage $f^{-1}\{z\}$, and $$\nu_1^f(z)=\nu_2^f(f(z))/\deg_zf.$$
By construction,  $ \f O_1^f\rightarrow \f O_2^f$
is a covering map between orbifolds.  
The following statement expresses the condition $g(\t S_f)\leq 1$ in terms 
of the Euler characteristic of $\f O_2^f$ (see \cite{gen}, Lemma 2.1).

\bl \l{ml} 
Let $f$ be a  rational function. Then $g(\t S_f)=0$ if and only if  $\chi(\f O_2^f)> 0$, and  $g(\t S_f)=1$ if and only if  $\chi(\f O_2^f)= 0$. \qed
\el

Using Lemma \ref{ml} one can show that Theorem \ref{main} is equivalent to the
following statement proved in the paper \cite{semi} (Theorem 6.1).


\bt \l{oip} Let $A,B,X$ be rational functions of degree at least two  such that $A\circ X=X\circ B$ and  $\C(B,X)=\C(z)$.   Then $\chi(\f O_1^X)\geq 0$, $\chi(\f O_2^X)\geq 0$, and 
the commutative diagram 
\be 
\begin{CD} \l{gooopa}
\f O_1^X @>B>> \f O_1^X\\
@VV X V @VV X V\\ 
\f O_2^X @>A >> \f O_2^X\ 
\end{CD}
\ee
consists of minimal holomorphic  maps between orbifolds. 
\et

Indeed,
a direct calculation shows that if $A,B,X$ is a primitive solution of \eqref{i1}, then  $A:\f O_1^X \rightarrow \f O_1^X$ and $B: \f O_2^X \rightarrow \f O_2^X$ are  minimal holomorphic  
maps between orbifolds (see \cite{semi}, Theorem 4.2). If 
Theorem \ref{main} is true, then Lemma \ref{ml} implies that $\chi(\f O_2^X)\geq 0$. Furthermore,
$\chi(\f O_1^X)\geq 0$, by \eqref{rh}.  
In turn, Theorem \ref{oip} implies Theorem \ref{main}, since  $\chi(\f O_2^X)\geq 0$ implies 
 $g(\t S_X)\leq 1$.  

\end{section}

\vskip 0.2cm
\noindent{\bf Acknowledgments}. The author is grateful to 
A. Eremenko  for discussions.

\end{document}